\documentclass[12pt,a4]{article}
\usepackage{amssymb}
\usepackage{amsfonts,amssymb}
\usepackage[]{graphicx}
\usepackage[cp1251]{inputenc}
\usepackage[english, russian]{babel}
\usepackage[margin=30mm]{geometry}

\begin{document}

\begin{flushleft}
 УДК 517.986.7:517.984.3
\end{flushleft}

\begin{center}
{\bf   А. Р. Миротин\\ О СОВМЕСТНЫХ СПЕКТРАХ НАБОРОВ
НЕОГРАНИЧЕННЫХ ОПЕРАТОРОВ}
\end{center}

\small{{\bf  Аннотация}.
Рассматривается несколько типов совместных спектров конечного семейства коммутирующих замкнутых операторов в банаховом пространстве, устанавливаются новые соотношения между этими спектрами (ранее было известно лишь включение  спектра Тейлора в коммутантный спектр), и для случая генераторов
 полугрупп доказываются  теоремы об отображении спектров, некоторые из которых являются обобщениями предшествующих результатов автора.  Даются  приложения полученных результатов к вопросу об  устойчивости многопараметрических полугрупп.

{\bf Ключевые слова}: неограниченный оператор, совместный спектр, функция Бернштейна, функциональное исчисление Бохнера-Филлипса,   теорема об отображении спектра, многопараметрическая полугруппа операторов, равномерная устойчивость, сильная устойчивость, абстрактная задача Коши.

\hspace{5mm}

{\bf  Abstract}. In this paper several joint spectra for a finite
commuting family of closed operators in Banach space are considered, some new relations between these spectra established (earlier only the inclusion of the Taylor spectrum in the commutant one was known), and
in the case of semigroup generators  spectral mapping theorems for such spectra are proved. Several of this theorems are generalizations of preceding results due to the author.  Applications to stability of multiparametric semigroups are given.
\hspace{5mm}

}

\begin{center}
\textbf{1. Введение и предварительные сведения}
\end{center}

 После основополагающих работ Р. Аренса, А. Кальдерона, Л. Валбрука и  Г. Е. Шилова  (см., например, \cite{GRS}, с. 91 ---  100)  различные виды совместных спектров  наборов ограниченных операторов   интенсивно изучались многими авторами (краткая история вопроса изложена в \cite{H}, гл. VI, см. также \cite{SZ} и библиографию там.)

     В данной работе рассмотрено  несколько типов совместных спектров  набора
  неограниченных  операторов в банаховом пространстве.   Установлены новые соотношения между этими спектрами. Показано, что в случае самосопряженных операторов в гильбертовом пространстве все рассматриваемые спектры (за исключением остаточного)  совпадают\footnote{Автор благодарит рецензента за постановку  задачи.}. Для случая
коммутирующих генераторов $C_0$-полугрупп  доказаны новые теоремы  о спектральных включениях и об отображении
этих спектров для многомерного функционального исчисления Бохнера-Филлипса, в том числе обобщающие некоторые результаты из \cite{SMZ} и \cite{FAN}. Даны приложения полученных результатов к вопросу об  устойчивости многопараметрических полугрупп и решений "каскада"\  абстрактных задач Коши.

 Для наборов неограниченных операторов совместный остаточный и совместный аппроксимативный спектры (равно как и их объединение) были введены в  \cite{SMZ}. Версии спектров Шилова,  Тейлора и коммутантного спектра,  рассматриваемые в данной работе, является  обобщениями соответствующих понятий из  \cite {Mir99} и \cite{FAN} соответственно, где рассматривались лишь генераторы полугрупп. Насколько известно автору, бикоммутантный спектр набора неограниченных  операторов  в банаховом пространстве ранее не рассматривался. Из нетривиальных соотношений между этими спектрами (см. теорему 1 ниже) ранее в \cite{FAN}  автором  было установлено лишь включение  спектра Тейлора в коммутантный спектр  для случая генераторов полугрупп.

Для однопараметрических полугрупп теорема об отображении спектра в исчислении Хилле-Филлипса опубликована в известной монографии этих авторов \cite{HiF}. По-видимому, первые теоремы об отображении спектра для одномерного исчисления Бохнера-Филлипса установлены  автором \cite{smz98}, результат следствия 7 данной статьи  был анонсирован, а также при дополнительных ограничениях доказан в \cite {SMZ},  теорема 4  усиливает теорему 3 из \cite {FAN}, первое равенство в утверждении 2)  теоремы 5 обобщает теорему 8.2 из \cite {Mir99}. Какие-либо аналоги или версии остальных теорем работы для случая наборов неограниченных операторов в банаховом пространстве автору не известны.

Отметим, что следствия 1 и 6 для однопараметрических полугрупп  известны давно (см. \cite{HiF}, \cite{OP}). Как справедливо было замечено рецензентом,  из более новых результатов в этом направлении следует указать слабые теоремы об отображении спектра для  групп. Первые результаты здесь (доказанные даже для представлений локально компактных абелевых групп, а потому верные и для $n$-параметрических групп операторов)  принадлежат А.~ Г.~ Баскакову \cite{Bas1}, \cite{Bas2}, по поводу  однопараметрических групп см. также \cite{EngNag}, \cite[глава 2]{Nee} (особенно теорему 2.4.4) \footnote{В последней работе излагается и история вопроса, однако пионерские работы А.~ Г.~ Баскакова не  упомянуты.}  и теорию Латушкина и Монтгомери-Смита \cite{Lat}.

Некоторые результаты об  устойчивости  многопараметрических полугрупп в общем контексте представлений подполугрупп локально компактных абелевых групп изложены в \cite{Batty}.

Таким образом, новизна статьи состоит в обобщении определений некоторых из введенных ранее совместных спектров, установлении  соотношений между рассматриваемыми совместными спектрами,  теорем об их отображении  и приложениях  к вопросам устойчивости многопараметрических полугрупп и решений задач Коши.  Обобщениями  доказанных ранее результатов автора являются теорема 4, часть второго утверждения теоремы 5, а также утверждение о включении спектра Тейлора в коммутантный в теореме 1.

Всюду далее $X$ есть комплексное
банахово пространство, $A_j, j=1,\dots,n$ --- коммутирующие замкнутые операторы в  $X$ с областями определения
$D(A_j)$,  резольвентные множества  $\rho(A_j)$ которых не пусты, $A=(A_1,\dots ,A_n)$;   коммутирование операторов  понимается как коммутирование  семейств  резольвент $R(\lambda_1,A_1),\dots,$ $
R(\lambda_n,A_n)$,  $\lambda_j\in \rho(A_j)$. Через $T(u)=T_1(u_1)\dots
T_n(u_n)$ $(u\in \Bbb{R}_+^n)$ будем обозначать  ограниченную
$n$-параметрическую $C_0$-полугруппу операторов  в пространстве $X$,
 $T_1, \dots ,T_n$ --- попарно коммутирующие
ограниченные однопараметрические $C_0$-полугруппы  в  $X,\ \|T_j(t)\|\leq M$.  Генератор полугруппы $T_j$ (частный генератор полугруппы $T$)  будет обозначаться
$A_j$. В случае генераторов коммутирование операторов $A_1,\dots,
A_n$  равносильно коммутированию операторов  $T_j(t_j); t_j>0, j=1,\dots,n$  (доказательство этого факта аналогично доказательству теоремы 16.2.1 из \cite{HiF}). Через
${\rm Gen}(X)$ мы будем обозначать множество  генераторов
 ограниченных $C_0$-полугрупп в $X$, а через $I$ --
единичный оператор. Отметим, что в случае генераторов векторное подпространство
$D(A):=\cap_{j=1}^nD(A_j)$ плотно в $X$.

Основы многомерного исчисления Бохнера-Филлипса  заложены
автором в    \cite {Mir99},  \cite{SMZ},     \cite {Mir97}  и \cite{Mir98}  (в \cite{SMZ} кратко изложена история вопроса и мотивировки).

 Напомним
необходимые понятия и факты из \cite {Mir97}  и \cite {Mir99}. Всюду ниже, где речь идет о многомерном исчислении Бохнера-Филлипса, $A=(A_1,\dots ,A_n)$ будет обозначать набор частных генераторов ограниченной $n$-параметрической $C_0$-полугруппы $T$ операторов  в пространстве $X$.

\textbf{ Определение 1.}   \cite {Mir97}  Будем говорить, что неположительная функция $\psi\in C^\infty((-\infty;0)^n)$ есть (отрицательная)
\textit{функция Бернштейна $n$ переменных} (или {\it принадлежит
классу ${\cal T}_n$}), если все ее частные производные первого
порядка абсолютно монотонны (т.е.  неотрицательны вместе со своими частными производными
всех порядков).

Ясно, что  ${\cal T}_n$ есть конус относительно поточечного
сложения функций и умножения на скаляр. Кроме того, в силу теоремы
6.1 из \cite {Mir99},  $\cup_{n=1}^\infty {\cal T}_n$ есть операда.

Известно \cite {Mir99}, \cite {MZ2013}, что каждая функция $\psi\in {\cal T}_n$ допускает
интегральное представление

\begin{equation}\label{psiots}
\psi(s)=c_0+c_1\cdot s+\int\limits_{\Bbb{R}_+^n\setminus \{0\}}
(e^{s\cdot u}-1)\,d\mu (u)   (s\in(-\infty;0)^n),
\end{equation}

\noindent где $c_0=\psi(-0):=\lim\limits_{s\to -0}\psi(s)$ (запись $s\to -0$ означает,
что $s_1\to -0, \dots, s_n\to -0$)
, а
$c_1$ из $\Bbb{R}_+^n$ и положительная мера $\mu$ на
$\Bbb{R}_+^n\setminus \{0\}$ определяются по $\psi$ однозначно (здесь и ниже $c\cdot s=\sum_{j=1}^nc_1^js_j$ при $c, s\in \mathbb{C}^n$).
Кроме того, $\psi$  продолжается по формуле (1) до функции,
голоморфной в области $\{\mathrm{Re}s<0\}\subset \mathbb{C}^n$  и
непрерывной в замыкании этой области.

\textbf{ Определение  2.}   \cite {Mir99}  Определим значение
функции $\psi$ из ${\cal T}_n$ вида (1) на наборе $A=(A_1,\dots
,A_n)$ при $x\in D(A)$ формулой
\begin{equation}\label{psiota}
\psi(A)x=c_0x+c_1\cdot Ax+\int\limits_{\Bbb{R}_+^n\setminus \{0\}}
(T(u)-I)x\,d\mu(u),
\end{equation}
где $c_1\cdot
Ax:=\sum_{j=1}^nc_1^jA_jx$ (интеграл понимается в смысле Бохнера).

Пусть $\psi\in{\cal T}_n, t\geq 0$. Тогда функция
$g_t(s):=e^{t\psi(s)}$ будет абсолютно монотонной на
$(-\infty;0)^n$, и $g_t(s)\leq 1.$ В силу многомерного варианта
теоремы Бернштейна-Уиддера  существует такая единственная
ограниченная положительная  мера $\nu_t$ на $\Bbb{R}_+^n,$ что при
$s\in (-\infty;0)^n$
\begin{equation}\label{gt}
g_t(s)=\int\limits_{\Bbb{R}_+^n} e^{s\cdot u}\,d\nu_t(u)={\cal
L}\nu_t(-s)
\end{equation}
 (здесь и ниже ${\cal L}\nu$
обозначает $n$-мерное преобразование Лапласа меры  $\nu$).

 \textbf{Определение 3.}  Используя обозначения,
введенные выше, положим ($x\in X$)
\begin{equation}\label{gt(A)}
g_t(A)x=\int\limits_{\Bbb{R}_+^n} T(u)x\,d\nu_t(u)
\end{equation}
(интеграл  Бохнера).

Отображение  $g(A):t\mapsto g_t(A)$ есть  ограниченная $C_0$-полугруппа
операторов в $X$, генератор которой есть замыкание оператора  $\psi(A)$, определенного выше лишь на $D(A)$. Это замыкание также обозначается $\psi(A)$.  В
одномерном случае $g(A)$  называется {\it полугруппой, подчиненной
полугруппе $T$} (терминология восходит к теории случайных процессов; по поводу одномерного случая см., например, \cite{App}, по поводу многомерного --- \cite{MZ2013}).

\begin{center}
\textbf{2. Совместные  спектры и соотношения между ними}
\end{center}

Далее мы будем иметь дело со следующими совместными спектрами
набора $A=(A_1,\dots,A_n)$ (вообще говоря, неограниченных)
 замкнутых операторов в комплексном банаховом пространстве $X$.

\textbf{ Определение  4. }  \cite{SMZ} {\it Совместный
аппроксимативный спектр} $\sigma_a(A)$ набора $A$   есть множество
тех $\lambda\in \Bbb{C}^n$, для которых найдется такая
последовательность векторов  $x_m\in D(A), \|x_m\|=1$, что $\sum_{j=1}^n\|A_jx_m-\lambda_jx_m\|\to 0\ (m\to\infty)$.

\textbf{ Определение  5. }  \cite{SMZ} {\it Совместный остаточный
спектр} $\sigma_R(A)$ набора $A$   есть множество тех $\lambda\in
\Bbb{C}^n$, для которых векторное пространство $\sum_{j=1}^n
\mathrm{Im}(\lambda_j I-A_j)$ не плотно в $X$.

 \textbf{Определение  6.}  \cite{SMZ} {\it Совместный  спектр}
$\sigma_J(A)$ набора $A$    определим равенством
\[
\sigma_J(A)=\sigma_a(A)\cup\sigma_R(A).
\]

Для определения еще трех  совместных спектров набора $A$   операторов в пространстве $X$
обозначим (мы предполагаем, что все резольвентные множества $\rho(A_j);j=1,\ldots ,n$ не пусты) через $\mathfrak{C}$ и  $\mathfrak{B}$  коммутант  (соответственно бикоммутант)
 семейства операторов ${\frak R}:=\{R(\lambda_j,A_j): \lambda_j\in \rho(A_j);j=1,\ldots ,n\}$
 в алгебре ${\frak L}(X)$ всех
ограниченных операторов в $X$. Алгебра ${\frak B}$ есть слабо (а потому и сильно)
замкнутая коммутативная наполненная подалгебра алгебры ${\frak
L}(X)$ с единицей, содержащая все операторы семейства ${\frak R}$. В случае, когда $A_j$ являются генераторами попарно коммутирующих полугрупп $T_j, j=1,\ldots ,n$ класса $C_0$ в пространстве $X$,  ${\frak B}$  содержит все операторы семейства $\{T_j(t): t\geq 0, j=1,\ldots ,n\}$  и совпадает с бикоммутантом этого семейства (доказательство повторяет доказательство теоремы 16.2.1 из \cite{HiF}).

Следующие два типа совместных спектров, наряду с совместным
аппроксимативным спектром,  для случая ограниченных операторов рассматривались в \cite{SZ}.

\textbf{ Определение  7.}  \cite{FAN} Будем говорить, что точка
$\lambda = (\lambda_1, \ldots, \lambda_n)\in \Bbb{C}^n$
 принадлежит {\it коммутантному спектру} набора $A$  (пишем
$\lambda\in\sigma'(A))$, если и только если не существует таких операторов
 $B_j\in \mathfrak{C}$, что

 \begin{equation}
 \sum_{j=1}^n(\lambda_jI - A_j)B_j = I.
 \end{equation}

 \textbf{Определение  8. }  Скажем, что точка
$\lambda \in \Bbb{C}^n$
принадлежит {\it бикоммутантному спектру} набора $A$  (пишем
$\lambda\in\sigma^{\prime\prime}(A))$,
 если  и только если равенство  (5) не может выполняться с операторами $B_j\in\mathfrak{B}$.

 Для определения еще одного из   рассматриваемых в данной работе совместных спектров как и в
\cite [теорема 16.3.1]{HiF}  получаем c помощью теорем 5.8.4 и 5.8.5 той же монографии, что при любом $j=1,\dots,n$ спектр Гельфанда ${\frak
M}_{\frak B}$ алгебры ${\frak B}$  разбивается на два подмножества
${\frak W}_j$ и ${\frak U}_j$, и существует непрерывная функция
$\alpha_j$ на ${\frak W}_j$ такая, что для любого комплексного
гомоморфизма $\varphi$ алгебры ${\frak B}$ имеем
\[
\varphi(R(z ,A_j))= \left\{\matrix {(z -\alpha_j(\varphi))^{-1}, &
\varphi\in{\frak W}_j, \cr 0, &\varphi\in{\frak U}_j \cr}\right.
\]
 для всех комплексных $z\in\rho(A_j)$, причем
$\sigma(A_j)= \alpha_j({\frak W}_j).$

\textbf{Определение  9.}  \cite {Mir99}  Используя предыдущие
обозначения, определим {\it совместный спектр Шилова} $\sigma(A)$ набора
$A$  замкнутых коммутирующих операторов в $X$   равенством
\[
\sigma(A)=\{(\alpha_1(\varphi),\ldots
,\alpha_n(\varphi)):\varphi\in \cap_{j=1}^n {\frak W}_j\}
\]
(мы считаем, что $\sigma(A)=\emptyset$, если
$\cap_{j=1}^n{\frak W}_j = \emptyset$).

\textbf{Лемма 1.}  1) \textit{ В случае одного оператора данные выше
определения совместных спектров совпадают с соответствующими
классическими.}

2) \textit{При $n>1$ спектры $\sigma_a(A)$ и
$\sigma_R(A)$ содержатся в
$\sigma_a(A_1)\times\dots\times\sigma_a(A_n)$  и
$\sigma_R(A_1)\times\dots\times\sigma_R(A_n)$, а  $\sigma_J(A)$,  $\sigma^{\prime\prime}(A),   \sigma(A)$  и
$\sigma'(A)$  --- в
$\sigma(A_1)\times\dots\times\sigma(A_n)$ соответственно}.

3)  \textit{Если все операторы $A_j$  ограничены, то  $\sigma(A)$   есть классический спектр Шилова набора  $A$ относительно  алгебры ${\frak B}$  и, стало быть,  совпадает с  $\sigma^{\prime\prime}(A)$.
}

Доказательство. Утверждения 1 и 2 очевидны. Если же все операторы $A_j$  ограничены, то в этом случае $A_j\in {\frak B}$, а потому $\varphi(R(z ,A_j))=1/(z-\varphi(A_j))$.  Следовательно, $\alpha_j(\varphi)=\varphi(A_j)$  при  $\varphi\in  {\frak W}_j$. Но  в нашем случае ${\frak U}_j=\emptyset$  (см.  \cite[теорема 16.3.1] {HiF}), а значит  $ {\frak W}_j={\frak M}_{\frak B}$. Это доказывает последнее утверждение леммы.

\textbf{Лемма 2. } 1) {\it Множества}  $\sigma_a(A), $
$\sigma'(A), \sigma^{\prime\prime}(A)$  \textit{и  $\sigma(A)$ замкнуты в $\mathbb{C}^n$}.

2)  \textit{Если операторы $A_j$ плотно определены, то  $\sigma(A)=\sigma(A')$, где $A':=(A_1',\dots,A_n')$ (штрих  обозначает
сопряженный оператор).
}

Доказательство. 1)  Пусть $\lambda\notin \sigma_a(A)$. Тогда существует такое
$\delta>0$, что при всех $x\in D(A),\ \|x\|=1$ справедливо
неравенство
\[
\sum\limits_{j=1}^n\|A_jx-\lambda_jx\|>\delta.
\]
Рассмотрим следующую окрестность точки $\lambda$:
\[
V=\{\lambda'\in
\mathbb{C}^n:|\lambda'_j-\lambda_j|<\frac{\delta}{2n},\
j=1,\dots,n\}.
\]
Для $\lambda'\in V$ и любых $x\in
D(A),\ \|x\|=1$  имеем
\[
\sum\limits_{j=1}^n\|A_jx-\lambda'_jx\|\geq
\sum\limits_{j=1}^n(\|A_jx-\lambda_jx\|-\|(\lambda_j-\lambda'_j)x\|)=
\]
\[
=\sum\limits_{j=1}^n\|A_jx-\lambda_jx\|-\sum\limits_{j=1}^n|\lambda_j-\lambda'_j|\geq
\delta-\frac{\delta}{2}=\frac{\delta}{2}.
\]
Следовательно, $\lambda'\notin \sigma_a(A)$.

 Докажем замкнутость коммутантного спектра.  Пусть $\lambda\notin\sigma'(A)$. Если точка $\lambda'$
принадлежит окрестности точки $\lambda$, задаваемой неравенством
\[
 \sum_{j=1}^n|\lambda_j - \lambda'_j|\|B_j\|<1,
\]
то оператор
\[
 C:=\sum_{j=1}^n(\lambda'_jI - A_j)B_j
 \]
принадлежит $\mathfrak{C}$ и имеет  обратный, принадлежащий  $\mathfrak{C}$,  поскольку оператор
\[
I- C=\sum_{j=1}^n(\lambda_j - \lambda'_j)B_j
\]
принадлежит $\mathfrak{C}$ и  $\|I-C\|<1$. Таким образом,
\[
\sum_{j=1}^n(\lambda'_jI - A_j)B_jC^{-1}=I,
\]
т. е. $\lambda'\notin\sigma'(A))$.
Замкнутость бикоммутантного спектра доказывается аналогично.

 Пусть теперь спектр $\sigma(A)$ не пуст, и $\lambda^0$ принадлежит
замыканию $\sigma(A)$, т.~ е. найдется такая последовательность
$\varphi_k\in \cap_{j=1}^n {\frak W}_j$, что при всех
$j=1,\dots,n$ имеем
$\alpha_j(\varphi_k))=z-\frac{1}{\varphi_k(R(z,A_j))}\to\lambda^0_j$ ($z\in\rho(A_j)$).
Переходя, если нужно, к подпоследовательности, можно считать, что
$\varphi_k\to\varphi_0\in{\frak M}_{\frak B}$. Тогда
$z-\frac{1}{\varphi_0(R(z,A_j))}=\lambda^0_j\ne\infty$, а потому
$\varphi_0\notin {\frak U}_j$ при всех $j$, и $\lambda^0_j=\alpha_j(\varphi_0)$ при всех $j$, т.
е. $\lambda^0\in\sigma(A)$, что и  завершает доказательство первого утверждения.

2) Обозначим через  $\mathfrak{B}^\prime$   бикоммутант
 семейства операторов $\{R(\lambda_j,A_j^\prime): \lambda_j\in \rho(A_j);j=1,\ldots ,n\}$
 в алгебре ${\frak L}(X^\prime)$ всех
ограниченных операторов в сопряженном пространстве $X^\prime$. Соотношения $\rho(A_j)=\rho(A_j'), R(\lambda_j,A_j)^\prime=R(\lambda_j,A_j^\prime)$, а также равенство   $(BC)'=C'B'\ (B, C\in {\frak L}(X))$, показывают, что отображение  $B\mapsto B'$ есть изоморфизм коммутативных банаховых алгебр  ${\frak B}$  и  ${\frak B}'$, при котором определение $\sigma(A)$ превращается в определение  $\sigma(A')$. Лемма доказана.

Напомним определение спектра Тейлора набора неограниченных операторов в пространстве $X$.

{\bf  Определение 10. \cite{FAN}} \textit{Пусть $\lambda\in\mathbb{C}^n$. Обобщенный комплекс Кошуля }${\cal
K}(X, \lambda I-A)$ набора $\lambda I-A$ определяется как последовательность пространств и
операторов
$$
0\longleftarrow
X_0\stackrel{d_0}{\longleftarrow}X_1\longleftarrow\cdots\stackrel{d_{n-1}}{\longleftarrow}X_n
\longleftarrow 0,
$$
 где $X_m=X\otimes\bigwedge^m \Bbb{C}^n\
(m=0,\ldots,n)$, а каждый оператор $ d_m=d_m(\lambda I-A)$ имеет область
определения $D(d_m)\subset X_{m+1}$, являющуюся линейной оболочкой
векторов вида $x\otimes e_{i_1}\wedge\cdots\wedge e_{i_{m+1}}$,
где $(e_k)_{k=1}^n$ -- стандартный базис в $\Bbb{C}^n$, $1\leq
i_1<\ldots < i_{m+1}\leq n$, а $x\in D(A_{j_1}\ldots
A_{j_{m+1}})$ для любой перестановки $(j_1,\ldots,j_{m+1})$
мультииндекса $(i_1,\ldots,i_{m+1})$. Для каждого $m=0,\ldots,
n-1$ оператор $d_m$ однозначно определяется на $D(d_m)$ равенством
$$
d_m(x\otimes e_{i_1}\wedge\cdots\wedge
e_{i_{m+1}})=\sum_{k=1}^{m+1}(-1)^k(\lambda_{i_k}I-A_{i_k})x\otimes
e_{i_1}\wedge\cdots\wedge\widehat{e_{i_k}}\wedge\cdots\wedge
e_{i_{m+1}}
$$
(здесь и ниже "крышка"\ , как обычно, означает, что стоящий под ней
сомножитель опущен).

Обобщенный комплекс Кошуля формально выглядит так же, как и
комплекс Кошуля для набора ограниченных операторов (см., например,
\cite{H}), но отличие состоит в том, что дифференциалы $d_m$ лишь замкнуты (но, вообще говоря, не
ограничены). При этом
 их области определения $D(d_m)$   выбраны таким
 образом, что

 $$
 d_m:D(d_m)\to D(d_{m-1}).
 $$

Легко проверить, что если набор $A$ замкнутых коммутирующих операторов в пространстве $X$ удовлетворяет следующему условию:
$$
\textit{\mbox{ если вектор } x \mbox{ принадлежит к }}
D(A_iA_j)\cap D(A_i),
$$
$$
\textit{\mbox{ то он принадлежит к }}
 D(A_jA_i), \mbox{ \textit{и} }
A_iA_jx=A_jA_ix\ (i,j=1,\ldots, n),\eqno{(K)}
$$
то последовательность   ${\cal K}(X,\lambda I-A)$ является комплексом в том
смысле, что  $d_{m-1}d_mx=0$  для всех $x\in D(d_m)$ и всех $m$.

Отметим, что набор частных генераторов многопараметрической полугруппы класса $(C_0)$  в пространстве $X$ удовлетворяет условию (К) в силу  \cite[теорема 10.9.4]{HiF}.

{\bf  Определение 11. \cite{FAN}} Пусть набор $A$ удовлетворяет  условию (К). \textit{Спектр Тейлора} $\sigma_T(A)$ этого набора
 состоит из тех $\lambda\in \Bbb{C}^n$, для которых комплекс ${\cal
K}(X,\lambda I-A)$ не является точным.

Для доказательства теоремы 1, описывающей соотношения между рассмотренными выше спектрами, нам потребуется следующая лемма, связывающая аппроксимативный спектр  со спектром Тейлора. Напомним, что \textit{совместный точечный спектр} $\sigma_p(A)$ набора $A$  есть множество тех $\lambda\in \mathbb{C}^n$, для которых найдется вектор $x\in  D(A),
x\ne 0$, удовлетворяющий равенствам $A_jx = \lambda_j x$ при всех $j = 1,\dots, n$.

\textbf{Лемма 3.} \textit{Соотношение $\lambda\in \sigma_a(A)$ выполняется тогда и только тогда, когда либо оператор   $d_{n-1}(\lambda I-A)$ не инъективен, либо его образ ${\rm Im}d_{n-1}(\lambda I-A)$ не замкнут.}

Доказательство. Пусть     $\lambda\in \sigma_a(A)$ и
последовательность векторов  $x_m\in D(A), \|x_m\|=1$ такова, что $\sum_{j=1}^n\|A_jx_m-\lambda_jx_m\|\to 0\ (m\to\infty)$.
Предположим, что  оператор   $d_{n-1}=d_{n-1}(\lambda I-A)$ инъективен, и докажем, что тогда  образ ${\rm Im}d_{n-1}$ не замкнут. Допустим противное, т. е. что ${\rm Im}d_{n-1}$ есть банахово подпространство пространства $X_{n-1}$. Далее для сокращения записи мы отождествим пространство $X_n$ с   $X$, а  пространство $X_{n-1}$ --- с прямой суммой $X\oplus\dots\oplus X$ посредством изоморфизмов $x\otimes e_{1}\wedge\cdots\wedge e_{n}\mapsto x$   и $\sum_{k=1}^n x_k\otimes e_{1}\wedge\cdots\wedge\widehat{e_{k}}\wedge\cdots\wedge
e_{n}\mapsto (x_k)_{k=1}^n$   соответственно. При этом отображение $d_{n-1}:D(A)\to X_{n-1}$ принимает вид $x\mapsto ((-1)^k(\lambda_k I-A_k)x)_{k=1}^n$. Обозначим через $X_A$ векторное пространство $D(A)$, наделенное нормой графика   $\|x\|_A:=\|x\|+\sum_{j=1}^n\|(\lambda_j I-A_j)x\|$. Поскольку оно банахово,  обратный оператор
$$
d_{n-1}^{-1}:{\rm Im}d_{n-1}\to X_A,\ ((-1)^k(\lambda_k I-A_k)x)_{k=1}^n\mapsto x
$$
ограничен по теореме о замкнутом графике. Следовательно, найдется такое  $C>0$, что при всех $x\in D(A)$ справедливо неравенство
$$
\|x\|+\sum\limits_{j=1}^n\|(\lambda_j I-A_j)x\|\leq C\sum\limits_{j=1}^n\|(\lambda_j I-A_j)x\|,
$$
Полагая здесь $x=x_m,\ m\to\infty$, приходим к противоречию.

Обратно, если оператор   $d_{n-1}$ не инъективен, то $\lambda\in\sigma_p(A)\subseteq\sigma_a(A)$. Пусть теперь этот оператор    инъективен, но его образ  не замкнут. Тогда обратный оператор $d_{n-1}^{-1}:{\rm Im}d_{n-1}\to X$ замкнут, но не ограничен ввиду незамкнутости своей области определения. Следовательно, для любого натурального  $m$ найдется такой вектор   $x_m\in D(A)$, что
\[
\|x_m\|>m\sum\limits_{j=1}^n\|(\lambda_j I-A_j)x_m\|
\]
   а потому    $\lambda\in \sigma_a(A)$. Лемма доказана.

\textbf{Теорема 1.}  \textit{1) Пусть набор $A$ замкнутых коммутирующих операторов в пространстве $X$ удовлетворяет  условию (К). Тогда справедливы включения}
\[
\sigma^{\prime\prime}(A)\supseteq \sigma(A)\supseteq \sigma'(A)\supseteq \sigma_T(A)\supseteq \sigma_J(A).
\]

\textit{2) Если исключить спектр Тейлора, то остальные включения в 1) останутся справедливыми для любого набора $A$ замкнутых коммутирующих операторов с непустыми резольвентными множествами.}

Доказательство. 1) Для обоснования первого включения   предположим, что $\sigma(A)\setminus\sigma^{\prime\prime}(A)\ne \emptyset.$ Тогда найдутся такие $B_j\in \mathfrak{B}$, что $\sum_{j=1}^n(\alpha_j(\varphi_0)I-A_j)B_j=I$ при некотором $\varphi_0\in\cap_{j=1}^n\mathfrak{W}_j$. Подставляя сюда $\alpha_j(\varphi_0)=z_j-1/\varphi_0(R(z_j,A_j))$, где $z\in\rho(A_j)$,
получим
\begin{equation}\label{sum}
\sum\limits_{j=1}^n\left(z_jB_j-A_jB_j-\frac{1}{\varphi_0(R(z_j,A_j))}B_j\right)=I.
\end{equation}
Операторы $A_jB_j$ определены на всем $X$, а потому ограничены по теореме о замкнутом графике. Следовательно, они  принадлежат $\mathfrak{B}$
(поскольку  $\mathfrak{B}$ принадлежат операторы $R(z_j,A_j)$, а с ними и  $R(z_j,A_j)B_j$).
Применяя теперь комплексный гомоморфизм $\varphi_0$ алгебры $\mathfrak{B}$ к обеим частям (\ref{sum}), получим $0=1$, и первое включение доказано.

Для доказательства второго включения допустим, что можно выбрать $\lambda\in\sigma'(A)\setminus\sigma(A)$. Зафиксируем  $\mu\in \rho(A_1)\times\dots\times\rho(A_n)$ и рассмотрим операторы $C_j:=I-(\mu_j-\lambda_j)R(\mu_j,A_j) \ (j=1,\dots, n)$. Ясно, что  $C_j\in\mathfrak{C}\cap\mathfrak{B}$, и что для любого комплексного гомоморфизма $\varphi$  алгебры   $\mathfrak{B}$ имеем, как в \cite[c. 204]{HiF},
\[
\varphi(C_j)=\left\{\begin{array}{rl}
\frac{\lambda_j-\alpha_j(\varphi)}{\mu_j-\alpha_j(\varphi)}\ne 0, & \mbox{ если } \varphi\in\mathfrak{W}_j\\
1, & \mbox{ если } \varphi\in\mathfrak{U}_j.
\end{array}\right.
\]
Следовательно, оператор $C_j$ имеет ограниченный обратный, который принадлежит алгебре  $\mathfrak{C}$ в силу наполненности последней. Кроме того, поскольку $C_j^{-1}$  коммутирует с резольвентами  $R(\zeta, A_k)$, этот оператор отображает векторное пространство $D(A_k)$  в себя $(j,k=1,\dots, n)$. Легко проверить также, что при всех $y\in D(A_j), \ j=1,\dots, n$ выполняется равенство $(\lambda_jI-A_j)y=(\mu_jI-A_j)C_jy$. Поскольку   $\mu\notin \sigma'(A)$, найдутся такие   $B'_j\in \mathfrak{C}$, что $\sum_{j=1}^n(\mu_j I-A_j)B'_jx=x$  при всех   $x\in X$. Если мы  положим теперь   $B_j:=C_j^{-1}B'_j$, то   $B_j\in \mathfrak{C}$, и     $\sum_{j=1}^n(\lambda_j I-A_j)B_j=I$, что противоречит выбору   $\lambda$.

 Для доказательства соотношения $\sigma'(A)\supseteq \sigma_T(A)$ заметим сначала, что если ограниченный оператор $B:X\to X$ принадлежит алгебре $\mathfrak{C}$, то
  он отображает линейное многообразие $D(A_{i_1}\ldots A_{i_m})$ в себя
для любого набора индексов $1\leq i_1<\ldots < i_{m}\leq n$, и
$$
BA_{i_1}\ldots A_{i_m}x=A_{i_1}BA_{i_2}\ldots A_{i_m}x=\cdots=
A_{i_1}\ldots A_{i_m}Bx
$$
для любого $x\in D(A_{i_1}\ldots A_{i_m})$. В самом деле, так как $B$ коммутирует со всеми операторами $A_j$, то указанное свойство справедливо при $m=1$ и легко доказывается индукцией по $m$. Таким образом, для набора $A$ справедливо заключение леммы 2 из \cite{FAN}. Но тогда для него справедлива и лемма 3 указанной работы, так как ее доказательство опирается лишь на условие (К) и лемму 2. Осталось заметить, что доказательство соотношения   $\sigma'(A)\supseteq \sigma_T(A)$, составляющее содержание теоремы 4 из  статьи \cite{FAN} (спектр   $\sigma'(A)$ обозначен там через   $\sigma_{\mbox {ш}}(A)$), опиралось лишь на леммы 2 и 3 этой статьи.

Перейдем к доказательству последнего утверждения части 1). Включение $\sigma_T(A)\supseteq \sigma_R(A)$ прямо
следует из определений, поскольку точность комплекса Кошуля в нулевом члене равносильна равенству $\sum_j{\rm Im}(\lambda_j I-A_j)=X$.
Осталось доказать, что $\sigma_T(A)\supseteq \sigma_a(A)$. Но при     $\lambda\in \sigma_a(A)$ из леммы 3 следует, что если   комплекс ${\cal K}(X,\lambda I-A)$ точен в члене   $X_n$, т. е. если оператор   $d_{n-1}$ инъективен, то  он не точен в члене   $X_{n-1}$, поскольку образ ${\rm Im}d_{n-1}$ не замкнут (ядро замкнутого оператора   замкнуто).

2) Поскольку доказательство первых двух включений не использовало свойство (К), осталось показать, что $\sigma_J(A)\subseteq \sigma'(A)$. Включение $\sigma_R (A)\subseteq \sigma'(A)$ следует из того, что для любого  $\lambda\notin\sigma'(A)$
$$\sum_{j=1}^n
\mathrm{Im}(\lambda_j I-A_j)\supseteq\left\{\sum_{j=1}^n
(\lambda_j I-A_j)B_jx:x\in X\right\}=X.
$$

 Если теперь предположить, что существует $\lambda\in\sigma_a(A)\setminus\sigma'(A)$, и  $(x_n)\subset D(A)$ --- соответствующий $\lambda$ совместный аппроксимативный собственный вектор, то, устремляя $n$ к бесконечности в равенстве   $\sum_{j=1}^n
B_j(\lambda_j I-A_j)x_n=x_n\ (B_j\in \mathfrak{C})$, получаем противоречие.
Теорема доказана.

\textbf{Замечание 1.}  Все фигурирующие в теореме 1  включения, за исключением первого, могут быть строгими даже в случае ограниченных операторов. В самом деле, в  \cite{T70}  построен пример, в котором $\sigma'(A)\ne \sigma_T(A)$.
С использованием этого примера в \cite[c. 144]{SZ} показано, что неравенство $\sigma'(A)\ne \sigma^{\prime\prime}(A)$ также возможно. Возможность неравенства $\sigma_J(A)\ne \sigma_T(A)$ (даже для ограниченных операторов в гильбертовом пространстве) установлена в следствии 3.4 работы \cite{Slod}, и осталось заметить, что в ограниченном случае $\sigma^{\prime\prime}(A)=\sigma(A)$.

 Ниже будут даны достаточные условия равенства $\sigma^{\prime\prime}(A)= \sigma(A)$ для неограниченных операторов. Вместе с тем, трудно усомниться, что первое включение в теореме 1 может быть строгим, но примером автор не располагает.

\textbf{Замечание 2.}  Пусть $\sigma_*$ обозначает один из семи рассматривавшихся выше совместных спектров, и пусть $\mathfrak{G}$  есть некоторый класс замкнутых операторов в пространстве $X$.  Будем говорить, что спектр $\sigma_*$ \textit{обладает проекционным свойством в классе } $\mathfrak{G}$, если $\sigma_*(A_1,\dots,A_{n})$  есть проекция $\sigma_*(A_1,\dots,A_{n+1})$ на первые $n$  координат в $\mathbb{C}^{n+1}$ для любых перестановочных $A_1,\dots,A_{n+1}\in\mathfrak{G}$ и  любого  $n$.   Ни один из рассмотренных в данной работе совместных спектров не обладает проекционным свойством в классе  генераторов ограниченных однопараметрических $C_0$-полугрупп, даже если ограничиться полугруппами, становящимися непрерывными в равномерной операторной  топологии.  В cамом  деле, достаточно взять $T(u_1,u_2)=T_1(u_1)$, где  $T_1$ --- нильпотентная полугруппа (относительно последней см., например, \cite[глава 1,  пример 9.6]{Goldst}).  В этом случае $\sigma_*(A_1,A_{2})\subset\sigma(A_1)\times\sigma(A_2)=\emptyset$, так как  $\sigma(A_1)=\emptyset$, и в то же время $\sigma(A_2)\ne\emptyset$. В  \cite{SZ} показано, что коммутантный и бикоммутантный спектры  не обладают проекционным свойством даже в классе ограниченных операторов, а аппроксимативный спектр в этом классе проекционным свойством обладает. Вопрос о наличии проекционного свойства у совместного аппроксимативного спектра в классе генераторов  однопараметрических $C_0$-полугрупп,   непрерывных  на  $(0,\infty)$ в равномерной операторной  топологии, насколько известно автору, остается открытым.

\begin{center}
\textbf{3. Случай самосопряженных операторов}
\end{center}

\textbf{  Теорема 2.} \textit{Для набора $A$ коммутирующих самосопряженных  операторов в комплексном гильбертовом пространстве $H$ справедливы следующие соотношения:}
\[
\sigma^{\prime\prime}(A)= \sigma(A)= \sigma'(A)=\sigma_T(A)= \sigma_J(A)=\sigma_a(A)\supseteq \sigma_R(A)=\sigma_p(A).
\]

Доказательство. Операторы набора $iA:=(iA_1,\dots,iA_n)$ являются частными генераторами унитарной группы $U(t):=\exp(it_1A_1)\dots\exp(it_nA_n)\ (t=(t_1,\dots,t_n))$, а потому удовлетворяют условию (К) \cite[теорема 10.9.4]{HiF}. Следовательно, этому условию удовлетворяют и операторы набора $A$  (другое доказательств этого факта вытекает из представления операторов $A_j$  операторами умножения на функции, полученного ниже).   Таким образом, ввиду теоремы 1 для доказательства справедливости первых пяти равенств достаточно убедиться, что $\sigma^{\prime\prime}(A)\subseteq \sigma_a(A)$.  Пусть $\lambda\in \mathbb{R}^n$.  Рассмотрим оператор
$$
S=\sum\limits_{j=1}^n(\lambda_jI-A_j)^2
$$
и покажем, что он самосопряжен. В самом деле, операторы $(\lambda_j I-A_j)^2$ самосопряжены и неотрицательны. Поэтому операторы $I+(\lambda_jI-A_j)^2$  имеют ограниченные обратные $C_j$, которые самосопряжены и коммутируют. Следовательно, самосопряжен и оператор
$$
S=\sum\limits_{j=1}^nC_j^{-1}-nI=\left(\sum\limits_{j=1}^nC_1\dots\widehat{C_j}\dots C_n\right)(C_1\dots C_{n})^{-1}-nI.
$$

Далее, если $\lambda\in \sigma^{\prime\prime}(A)$, то  $0\in \sigma(S)$, поскольку в противном случае мы имели бы равенство
$$
\sum\limits_{j=1}^n(\lambda_j I-A_j)((\lambda_j I-A_j)S^{-1})=I,
$$
противоречащее выбору $\lambda$ (операторы  $B_j:=(\lambda_j I-A_j)S^{-1}$ ограничены по теореме о замкнутом графике и, стало быть, принадлежат $\mathfrak{B}$ вслед за $S^{-1}$).

Так как $\sigma_a(S)=\sigma(S)$ (оператор $S$ самосопряжен),  найдется такая
последовательность векторов  $x_m\in D(S), \|x_m\|=1$, что $\|Sx_m\|\to 0\ (m\to\infty)$.
Для завершения доказательства первых пяти равенств теперь нужно лишь заметить, что
$$
\sum\limits_{j=1}^n\|(\lambda_j I-A_j)x_m\|^2=\left\langle\sum\limits_{j=1}^n(\lambda_jI-A_j)^2x_m,x_m\right\rangle\leq \|Sx_m\|
$$
(угловые скобки обозначают скалярное произведение в $H$), а потому $\lambda\in \sigma_a(A)$.

Поскольку $\sigma_p(A)\subseteq \sigma_a(A)$,  осталось установить лишь последнее равенство. Но включение $\lambda\in \sigma_R(A)$ равносильно тому, что найдется такой ненулевой вектор   $x\in H$, что    $\langle\sum_{k=1}^n(\lambda_kI-A_k)y_k,x\rangle=0$ при всех    $y_k\in D(A_k)$, что, в свою очередь,   равносильно тому, что
$(\lambda_kI-A_k)x=0$ при всех    $k=1,\dots,n$. Теорема доказана.

 \textbf{Замечание 3.}  Для случая ограниченных коммутирующих  нормальных операторов равенство  $\sigma^{\prime\prime}(A)$ $=\sigma_a(A)$  было установлено в \cite{Dash}.

Для того чтобы получить еще одну характеризацию совместного спектра набора коммутирующих самосопряженных операторов, нужна некоторая подготовка. Ниже в этом разделе $H$ есть сепарабельное гильбертово пространство. Поскольку $\mathfrak{B}$ есть коммутативная алгебра фон Неймана \cite{Dix}, существует упорядоченное представление пространства $H$ относительно $\mathfrak{B}$ (см., например, \cite[c. 17 -- 18]{Schw}). А именно, найдется такая убывающая последовательность борелевских множеств $e_k\subset \mathbb{R}, e_1= \mathbb{R}$ и такая конечная регулярная борелевская мера $\mu$ на $\mathbb{R}$, что $H$ изоморфно гильбертовой прямой  сумме $\widetilde H=\bigoplus\{L^2(e_k,\mu):k\in \mathbb{N}\}$, и  этот изоморфизм индуцирует изоморфизм $V$  алгебры $\mathfrak{B}$ на подалгебру $\widetilde{\mathfrak{B}}$ алгебры  $\mathfrak{L}(\widetilde H)$, состоящую из операторов вида  $M(b):(f_k)_{k=1}^\infty\mapsto (bf_k)_{k=1}^\infty$, где  $f_k\in L^2(e_k,\mu)$, а $b$ пробегает алгебру $B(\mathbb{R})$   ограниченных борелевских функций на  $\mathbb{R}$. Пусть $r_j^t\in B(\mathbb{R})$ таковы, что оператор $M(r_j^t)$ соответствует при этом изоморфизме оператору $R(t,A_j)\ (t\in\rho(A_j); j=1,\dots,n)$. Поскольку оператор
$M(r_j^t)$  самосопряжен и инъективен  вслед за $R(t,A_j)$,  функция $r_j^t$ вещественнозначна и отлична от нуля $\mu$-~п.~в.  Кроме того, первое резольвентное тождество влечет равенство $r_j^t-r_j^s=(s-t)r_j^tr_j^s\ (t,s\in\rho(A_j); j=1,\dots,n)$. Отсюда следует, что, если мы для каждого $j=1,\dots,n$ и  $t\in\rho(A_j)$ положим $a_j:=t-1/r_j^t$, то эта вещественнозначная борелевская функция не зависит от  $t$. Введем в рассмотрение операторы  $\widetilde A_j$ в   $\widetilde H$, соответствующие операторам $A_j$, следующим образом:
$$
D(\widetilde A_j):=\{f=(f_k)_{k=1}^\infty\in \widetilde H: a_jf_k\in  L^2(e_k,\mu) \mbox{ при }  k\in \mathbb{N}\},\ A_jf:=(a_jf_k)_{k=1}^\infty \ (j=1,\dots,n).
$$
Эти операторы самосопряжены, и $\rho(\widetilde A_j)=\rho(A_j)$. В самом деле, хорошо известно, что сужения $\widetilde A_j|L^2(e_k,\mu)$ являются самосопряженными операторами. Далее, при $t\in \rho(A_j)$ легко проверить, что оператор $tI-\widetilde A_j$  обратим и $R(t,\widetilde A_j)=M(r_j^t)$, а потому  $\rho(\widetilde A_j)\supseteq\rho(A_j)$, и при изоморфизме $V$ оператор  $R(t, A_j)$ переходит  в  $R(t, \widetilde A_j)$, а значит $\sigma(R(t, A_j))=\sigma(R(t, \widetilde A_j))$. Применяя теперь теорему об отображении спектров для резольвенты, получаем
$$
\left\{\frac{1}{t-\xi}:\xi\in\sigma(A_j)\right\}=\left\{\frac{1}{t-\zeta}:\zeta\in\sigma(\widetilde A_j)\right\},
$$
откуда следует, что $\sigma(A_j)=\sigma(\widetilde A_j)$. Если мы теперь положим  $\widetilde A:=(\widetilde A_1,\dots,\widetilde A_n)$, то для спектров Шилова будем иметь  $\sigma(A)=\sigma(\widetilde A)$, поскольку, как показано выше, эти спектры   определяются в терминах, инвариантных относительно $V$. (Из описанной выше конструкции также вытекает, что набор самосопряженных коммутирующих операторов удовлетворяют условию (К).)

Для получения требуемой характеризации, введем теперь, следуя \cite{Dash},

\textbf{Определение 12. }  \textit{Совместной существенной областью значений} набора функций $a=(a_1,\dots,a_n)$    называется множество $\mathcal{E}(a)$ таких $\beta\in \mathbb{R}^n$, что при любом $\varepsilon>0$  множество $\{s\in \mathbb{R}:\sum_{j=1}^n|\beta_j-a_j(s)|<\varepsilon\}$  имеет положительную $\mu$-меру.

\textbf{  Теорема 3.} \textit{Во введенных выше обозначениях справедливо равенство} $\sigma(A)=\mathcal{E}(a)$.

Доказательство. В силу теоремы 2 и равенства $\sigma(A)=\sigma(\widetilde A)$ для доказательства включения $\sigma(A)\subseteq\mathcal{E}(a)$   достаточно показать, что   $\sigma_a(\widetilde A)\subseteq\mathcal{E}(a)$. Допустим противное и выберем  $\lambda\in \sigma_a(\widetilde A)\setminus\mathcal{E}(a)$. Тогда найдется такая последовательность   $f_m\in D(\widetilde A)$, что   $\|f_m\|_{\widetilde H}=1$, но   $\|(\lambda_j I-\widetilde A_j)f_m\|_{\widetilde H}\to 0\ (m\to\infty; j=1,\dots,n)$.  С другой стороны, найдется такое   $\varepsilon_0>0$, что   множество  $\{s\in \mathbb{R}:\sum_{j=1}^n|\lambda_j-a_j(s)|<\varepsilon_0\}$   имеет     меру нуль, т. е. $\sum_{j=1}^n|\lambda_j-a_j(s)|\geq\varepsilon_0$  $\mu$~-п.~в., а потому
$$
\sum\limits_{j=1}^n|\lambda_j-a_j(s)|^2\geq\frac{1}{n}\left(\sum\limits_{j=1}^n|\lambda_j-a_j(s)|\right)^2\geq\frac{\varepsilon_0^2}{n} \  \mu-\mbox{п. в.}
$$

Следовательно,
$$
\sum\limits_{j=1}^n \|(\lambda_j I-\widetilde A_j)f_m\|_{\widetilde H}^2=\sum\limits_{j=1}^n \sum\limits_{k=1}^\infty \|(\lambda_j-a_j)f_{mk}\|_{L^2(e_k,\mu)}^2=\sum\limits_{k=1}^\infty\sum\limits_{j=1}^n
\int\limits_{e_k}|\lambda_j-a_j(s)|^2|f_{mk}(s)|^2d\mu(s)=
$$
$$
=\sum\limits_{k=1}^\infty\int\limits_{e_k}\left(\sum\limits_{j=1}^n
|\lambda_j-a_j(s)|^2\right)|f_{mk}(s)|^2d\mu(s)\geq \frac{\varepsilon_0^2}{n}\sum\limits_{k=1}^\infty\|f_{mk}\|_{L^2(e_k,\mu)}^2=\frac{\varepsilon_0^2}{n},
$$
и мы получили противоречие.

Для доказательства обратного включения, снова в силу теоремы 2, достаточно показать, что   $\mathcal{E}(a)\subseteq\sigma^{\prime\prime}(\widetilde A)$. Допустим противное и выберем  $\lambda\in \mathcal{E}(a)\setminus\sigma^{\prime\prime}(\widetilde A)$. Тогда найдутся такие $b_j\in B(\mathbb{R})$, что $\sum_{j=1}^n (\lambda_j I-\widetilde A_j)M(b_j)=I$, а потому $\sum_{j=1}^n (\lambda_j -a_j(s))b_j(s)=1$ $\mu$-~п.~ в. Пусть константа $C>0$ такова, что $|b_j(s)|\leq C$ при всех  $s\in \mathbb{R}, j=1,\dots,n$.
    Множество $E=\{s\in \mathbb{R}:\sum_{j=1}^n |\lambda_j -a_j(s)|<1/2C\}$   имеет положительную $\mu$-меру, и для $\mu$-п. в.  $s\in E$ имеем
$$
1=\left|\sum\limits_{j=1}^n (\lambda_j -a_j(s))b_j(s)\right|\leq C\sum\limits_{j=1}^n |\lambda_j -a_j(s)|<\frac{1}{2}.
$$
 Полученное противоречие завершает доказательство теоремы.

\textbf{ Замечания 4.} 1) В монографии  \cite{BS}  показано, что совместный аппроксимативный спектр  набора $A$ совпадает также с носителем произведения спектральных мер операторов этого набора. Таким образом, с учетом теорем 2 и 3 мы имеем восемь равносильных определений совместного спектра набора коммутирующих самосопряженных операторов в гильбертовом пространстве.

2) Для набора $A$ ограниченных коммутирующих  нормальных операторов  в \cite{Dash} доказано, что $\sigma^{\prime\prime}(A)=\mathcal{E}(a)$.

\begin{center}
\textbf{4. Теоремы об отображении совместных  спектров}
\end{center}

Как отмечено в \cite {SMZ}, уже для функции $\psi\in {\cal T}_1$ теорема
об отображении спектра, понимаемая как утверждение о справедливости равенства
$\psi(\sigma(A))=\sigma(\psi(A))$, где $A\in {\rm
Gen}(X),\sigma(A)$ -- спектр оператора $A$, может не выполняться.
Заметим, что нижеследующая теорема 4  о спектральном включении усиливает теорему 3 из \cite {FAN} (где утверждается лишь, что $\sigma(\psi(A))\supseteq \psi(\sigma^{\prime}(A))$) и доказывается аналогично.  Мы приводим доказательство для полноты изложения.

\textbf{Теорема  4.} \textit{ Для любой функции $\psi\in {\cal T}_n$  и набора $A$ попарно коммутирующих
генераторов ограниченных $C_0$-полугрупп в банаховом пространстве
$X$ имеет  место  включение}
\[
\sigma(\psi(A))\supseteq \psi(\sigma^{\prime\prime}(A)).
\]

Доказательство. Докажем сначала, что
\begin{equation}\label{rhocc}
\mathbb{C}^n\setminus\sigma^{\prime\prime}(A)\supseteq\{\lambda\in\Bbb{C}^n: {\rm
Re}\lambda<0, \psi(\lambda)\in\rho(\psi(A))\}.
\end{equation}
 При $\lambda\in\Bbb{C}^n, {\rm Re}\lambda<0$ и $x\in
D(A)$ имеем в силу (\ref{psiots}) и (\ref{psiota})
\begin{equation}\label{psix}
\psi(\lambda)x - \psi(A)x = c_1\cdot(\lambda I - A)x +
\int\limits_{\Bbb{R}_+^n} (e^{\lambda\cdot u}I -
T(u))xd\mu(u).
\end{equation}

Для любых ограниченных операторов $B_j, C_j, j=1,\ldots,n$ в $X$
справедливо тождество
$$
\prod_{j=1}^nB_j - \prod_{j=1}^nC_j = (B_1 - C_1)\prod_{j=2}^nB_j
+ C_1(B_2 - C_2)\prod_{j=3}^nB_j + \cdots +
\prod_{j=1}^{n-1}C_j(B_n - C_n).
$$
С учетом коммутирования полугрупп $T_j$ получаем отсюда, что
\begin{equation}\label{eliu}
e^{\lambda\cdot u}I-T(u)=\sum_{j=1}^n(e^{\lambda_ju_j}I-T_j(u_j))
U_j(u),
\end{equation}
 где
$$
U_j(u)=\prod\limits_{1\leq l<j}T_l(u_l)\prod\limits_{j<k\leq
n}e^{\lambda_k u_k}
$$
--- ограниченные операторы в $X$, коммутирующие с $T_i,
i=1,\ldots,n, \|U_j(u)\|\leq$  $M^{n-1}$.  Известно
(см., например, \cite{OP}, формулы (8.1a), (8.1b)) что
равенства
$$
V_j(u_j)x=\int\limits_0^{u_j}e^{(u_j-s)\lambda_j}T_j(s)xds
$$
определяют ограниченные операторы в $X, {\rm
Im(}V_j(u_j))\subseteq D(A_j),$ также коммутирующие со всеми
$T_i$, причем $\|V_j(u_j)x\|\leq M\|x\|(e^{u_j{\rm
Re}\lambda_j}-1)/{\rm Re}\lambda_j$ и
$$
(e^{\lambda_ju_j}I-T_j(u_j))x=(\lambda_jI-A_j)V_j(u_j)x \mbox{ при
} x\in X.
$$
 Подставляя это в (\ref{eliu}), а результат -- в (\ref{psix}), имеем
\[
(\psi(\lambda)I-\psi(A))x=\sum_{j=1}^nc_1^j(\lambda_jI-A_j)x+
\sum_{j=1}^n\int\limits_{\Bbb{R}_+^n}(\lambda_jI-A_j)V_j(u_j)U_j(u)xd\mu(u)=
\]
\begin{equation}\label{psiI}
\sum_{j=1}^n(\lambda_jI-A_j)W^\lambda_jx,
\end{equation}
 где операторы
$$
W^\lambda_jx=c_1^jx+\int\limits_{\Bbb{R}_+^n}V_j(u_j)U_j(u)xd\mu(u)
$$
\noindent ограничены  в $X$, так как
$$
\|W^\lambda_jx\|\leq
\left(c_1^j+M^{n-1}\int\limits_{\Bbb{R}_+^n}\|V_j(u_j)\|d\mu(u)\right)\|x\|\leq
$$
$$
\leq\left(c_1^j+\frac{M^n}{{\rm
Re}\lambda_j}\int\limits_{\Bbb{R}_+^n}(e^{u_j{\rm
Re}\lambda_j}-1)d\mu(u)\right)\|x\|,
$$
 (интегралы сходятся, поскольку последний интеграл выражается через
$\psi(({\rm Re}\lambda_j) e_j)$, где $(e_j)$ --- стандартный
базис в
 $\Bbb{R}^n$). Кроме того, ${\rm
Im(}W^\lambda_j)\subseteq D(A_j)$,
 и $W^\lambda_j$ коммутируют со всеми $T_i$, а потому и с $\psi(A)$.

 Если $\psi(\lambda)\in
\rho(\psi(A)),$ то, заменяя в (\ref{psiI}) $x$ на
$(\psi(\lambda)I-\psi(A))^{-1}x$,   получаем
$$
x=\sum_{j=1}^n(\lambda_jI-A_j)W^\lambda_j(\psi(\lambda)I-\psi(A))^{-1}x
\mbox { при } x\in X.
$$
 Заметим, что резольвента $(\psi(\lambda)I-\psi(A))^{-1}$  принадлежит $\mathfrak{B}$, поскольку этой алгебре принадлежит полугруппа $g(A)$, генератором которой является оператор  $\psi(A)$. Таким образом,  $\lambda\notin\sigma''(A),$ что  доказывает (\ref{rhocc}).

Рассмотрим теперь случай, когда $\sigma(A_j)\subseteq \{{\rm
Re}\lambda_j<0\}$ при всех $j$.  Если  мы предположим, что в этом
случае теорема неверна, т. е. найдется число из
$\psi(\sigma^{\prime\prime}(A)),$ не принадлежащее $\sigma(\psi(A)),$
то $\psi(\lambda)\in \rho(\psi(A))$ при некотором $\lambda\in
\sigma^{\prime\prime}(A),$ а это противоречит (\ref{rhocc}).

Перейдем к общему случаю. В силу ограниченности рассматриваемых
полугрупп  при всех $j=1,\ldots,n$ справедливы включения
$\sigma(A_j)\subseteq \{{\rm Re}\lambda_j\leq 0\}$. Тогда для
любого $\delta>0$ операторы $A_j-\delta I\in {\rm Gen}(X)$ и
$\sigma(A_j-\delta I)\subseteq \{{\rm Re}\lambda_j<0\}$ при всех
$j$. Кроме того, $\psi_\delta(s):=\psi(s-\bar\delta)\in {\cal
T}_n$ и $\psi_\delta(A)=\psi(A-\bar\delta I)$, где положено
$\bar\delta:=(\delta,\ldots, \delta)$. По доказанному выше
\begin{equation}\label{sigmapsi}
\sigma(\psi_\delta(A))=\sigma(\psi(A-\bar\delta I))\supseteq
\psi(\sigma^{\prime\prime}(A-\bar\delta
I))=\psi_\delta(\sigma^{\prime\prime}(A)).
\end{equation}
Далее,  в доказательстве  теоремы 4.1 из \cite{Mir99}  показано, что
$\psi(A)x=\psi_\delta(A)x+F_\delta x$ при $x\in D(A)$, где
$F_\delta$ -- ограниченный оператор. Но тогда
$\psi(A)=\psi_\delta(A)+F_\delta$, поскольку по теореме 4 из \cite{SMZ}
$D(A)$ есть существенная область для обеих частей этого равенства.
Кроме того,  в доказательстве  теоремы 4.1 из \cite{Mir99} установлено, что $\|F_\delta\|\to 0 \ (\delta\to
0)$, а потому по теореме о спектре возмущенного оператора (см.
\cite{Kato}, теорема IV.3.6; $\rho_H$ обозначает метрику Хаусдорфа)
$$
\rho_H(\sigma(\psi_\delta(A)),\sigma(\psi(A)))\to 0 \ (\delta\to
0).
$$
Следовательно, для любого  $\varepsilon>0$ при достаточно малых
$\delta>0$
$$
\sigma(\psi_\delta(A))\subseteq B_\varepsilon[\sigma(\psi(A))],
$$
где для множества $N\subset \Bbb{C}$ положено
$B_\varepsilon[N]:=\{z\in \Bbb{C}:\inf_{t\in N}|z-t|\leq
\varepsilon\}$. А тогда в силу (\ref{sigmapsi})
$$
\psi_\delta(\sigma^{\prime\prime}(A))\subseteq
B_\varepsilon[\sigma(\psi(A))].
$$

С другой стороны, при достаточно малых $\delta>0$
$$
\psi(\sigma^{\prime\prime}(A))\subseteq
B_\varepsilon[\psi_\delta(\sigma^{\prime\prime}(A))].
$$
Это следует из того, что
$\psi(\lambda)-\psi(\lambda-\bar\delta)\to 0\ (\delta\to +0)$
равномерно по $\lambda\in \{{\rm Re}\lambda\leq 0\}$, так как
$$
|\psi(\lambda)-\psi(\lambda-\bar\delta)|\leq
c_1\cdot\bar\delta+\int\limits_{\Bbb{R}_+^n}(1-e^{-u\cdot\bar\delta})d\mu(u)\to
0\ (\delta\to +0)
$$
по теореме Б. Леви.

Следовательно,
$$
\psi(\sigma^{\prime\prime}(A))\subseteq B_\varepsilon[
B_\varepsilon[\sigma(\psi(A))]]\subseteq
B_{2\varepsilon}[\sigma(\psi(A))],
$$
а потому
$$
\psi(\sigma^{\prime\prime}(A))\subseteq
\bigcap\limits_{\varepsilon>0}B_{2\varepsilon}[\sigma(\psi(A))]=\sigma(\psi(A))
$$
в силу замкнутости $\sigma(\psi(A))$. Теорема доказана.

При $\psi(s)=e^{u\cdot s}-1\ (u\in \mathbb{R}_+^n)$ получаем

{\bf Следствие 1.}  \textit{При  $u\in \mathbb{R}_+^n имеем$}
\[
\sigma(T(u))\supseteq e^{u\cdot\sigma^{\prime\prime}(A)}
\]
(\textit{здесь и ниже  для $\Omega\subset \mathbb{C}^n,\ u\in \mathbb{R}^n$ положено $u\cdot \Omega=\left\{\sum_{j=1}^n u_jz_j:z\in \Omega\right\}$}).

Идущее ниже следствие полезно сравнить со следствием 18  из  \cite{SMZ}.

{\bf Следствие 2.} \textit{Пусть оператор $\psi(A)$ ограничен хотя бы для одного набора $A$, удовлетворяющего условию  $\sigma^{\prime\prime}(A)\supseteq
(-\infty;a)^n$ с некоторым $a\leq 0$. Тогда функция $\psi$ ограничена на $(-\infty;0)^n$.}

 Доказательство. В  силу  теоремы 4   функция   $\psi$ ограничена на $(-\infty;a)^n$, т. е. при некотором $C<0$ для всех $r\in (-\infty;a)^n$ имеем $\psi(r)\geq C$. Для любого $s\in (-\infty;0)^n$  выберем $r\in (-\infty;a)^n$  так, что $r_j\leq s_j$  при всех $j$. Тогда,  поскольку  $\psi$ возрастает  по каждой переменной в отдельности, получим  $\psi(s)\geq\psi(r)\geq C$.

Первое равенство в утверждении 2) следующей теоремы об отображении спектров  обобщает теорему 8.2 из \cite {Mir99}, где рассмотрен лишь случай генераторов голоморфных полугрупп,  остальные утверждения являются новыми.

\textbf{  Теорема 5. }  1) {\it  Пусть $\psi\in{\cal
T}_n$, $A$ --- набор коммутирующих операторов из $\mathrm{Gen}(X)$.   Имеют место следующие
утверждения:

{\rm 1)} \centerline{$\sigma(\psi(A))\supseteq\psi(\sigma(A)).$}

 {\rm 2)}  Если все полугруппы $T_j$ и
$g(A)$ непрерывны на
$(0;+\infty)$ в равномерной операторной топологии,  операторы $A_j$ не ограничены и   $\psi(-\infty)=$ $-\infty$, то
\begin{equation}\label{sigmaotpsi}
\sigma(\psi(A))=\psi(\sigma(A))=\psi(\sigma^{\prime\prime}(A)).
\end{equation}

{\rm 3)}  Если все полугруппы $T_j$  непрерывны на
$(0;+\infty)$ в равномерной операторной топологии,  операторы $A_j$ не ограничены и $\psi(-\infty)\ne -\infty$, то}
\[
\sigma(\psi(A))=\psi(\sigma(A))\cup\{\psi(-\infty)\}.
\]

Доказательство.  Утверждение 1) вытекает из теорем 1 и 4.

Докажем утверждение 2).  Так как полугруппы $T_j$ и $g(A)$ непрерывны на $(0;+\infty)$
в равномерной операторной топологии, то  при всех $t>0$
\begin{equation}\label{etsigma}
e^{t\sigma(\psi(A))}=\sigma(g_t(A))\setminus\{0\}
\end{equation}
(см., например, \cite{EngNag} или \cite[предложение 8.5]{OP}).

 Для ограниченной регулярной борелевской меры $a$ на $\mathbb{R}_+^n$  и функции $\psi_1(s)=({\cal
L}a)(-s)$ справедливо
равенство
\begin{equation}\label{Psi}
\psi_1(\sigma(A))\cup \{a(\{0\})\}=\sigma(\Psi(a,A)),
\end{equation}
где  $\Psi(a,A):=\int_{{\Bbb R}_+^n}T(u)da(u)$  (см. \cite[раздел 8]{Mir99}).  Отсюда при
$a=\nu_t, \psi_1(s)=g_t(s),$ $ \Psi(a,A)=g_t(A)$ имеем
\begin{equation}\label{etpsi}
\sigma(g_t(A))=e^{t\psi(\sigma(A))}\cup\{0\},
\end{equation}
поскольку $\nu_t(\{0\})=g_t(-\infty)=0$ по тауберовой
теореме для преобразования Лапласа. Сопоставляя
равенства (\ref{etsigma}) и  (\ref{etpsi}), получаем, что при $t>0$
\begin{equation}\label{ee}
e^{t\sigma(\psi(A))}=e^{t\psi(\sigma(A))}.
\end{equation}
 Известно (см., например, \cite[с. 113, теорема 4.18] {EngNag}),
что множество  $\sigma(\psi(A))$, будучи спектром генератора $C_0$-полугруппы, которая  становится непрерывной в равномерной операторной топологии,  имеет ограниченное
пересечение с любой  правой полуплоскостью $\{b\leq s_1\}$
комплексной $s$-плоскости ($s=s_1+is_2$).
Фиксируем вещественное число $s_1^0$ и обозначим через $E$ и $P$
сечения множеств $\sigma(\psi(A))$ и $\psi(\sigma(A))$
соответственно  вертикальной прямой $\{s_1=s_1^0\}$.
Множество $E$ ограничено. Из (\ref{ee}) следует, что $e^{tiE}= e^{tiP}$.
Фиксируем $y\in P$.  В силу  последнего равенства,  для
любого $t>0$ найдется такой элемент $x_t\in E$ и целое число
$n_t$, что $t(y-x_t)=2\pi n_t.$ Поскольку левая часть этого
равенства есть бесконечно малая при $t\to 0$, то $n_t=0$ при
достаточно малых $t$, а потому $y=x_t\in E$ при этих $t$. Значит,
$P\subseteq E$. В частности, множество $P$ тоже ограничено. Но тогда
аналогично получаем, что $E\subseteq P$, т. е. $E=P$.  В силу
произвольности $s_1^0$, отсюда следует первое из  равенств (\ref{sigmaotpsi}).  Наконец,  применяя последовательно  это
равенство,  теорему 4 и теорему 1,  имеем
\[
\psi(\sigma(A))=\sigma(\psi(A))\supseteq \psi(\sigma^{\prime\prime}(A))\supseteq\psi(\sigma(A)),
\]
что и завершает доказательство утверждения 2).

Докажем 3).   По теореме 8.1 из \cite{Mir99}  $\psi(s)=({\cal L}a)(-s)$, где $a$  ---  ограниченная мера на $\mathbb{R}_+^n$, причем $\psi(A)=\Psi(a,A)$. Тогда в силу (\ref{Psi}) $\sigma(\psi(A))=\psi(\sigma(A))\cup a(\{0\})$, и  осталось заметить, что по тауберовой теореме для преобразования Лапласа $a(\{0\})=\psi(-\infty)$.

Покажем, что утверждение 3) предыдущей теоремы может быть неверным для случая ограниченных генераторов.

\textbf{Пример 2. } Положим $n=1, \psi(s_1)=e^{s_1}-1$, и пусть $A_1$ есть ограниченный оператор в $X$. В этом случае множество $\sigma(\psi(A_1))=\sigma(e^{A_1})-1$ не содержит числа $\psi(-\infty)=-1$, поскольку оператор $e^{A_1}$ ограниченно обратим, и, стало быть,  утверждение 3) неверно.

\textbf{Следствие 3.}  \textit{Если все полугруппы $T_j$  непрерывны на
$(0;+\infty)$ в равномерной операторной  топологии, а их генераторы $A_j$ не ограничены, то для любого $c_1\in \mathbb{R}_+^n$   справедливы равенства
\[
\sigma\left(\overline{\sum\limits_{j=1}^n c_1^j A_j}\right)=c_1\cdot\sigma(A)=c_1\cdot\sigma^{\prime\prime}(A),
\]
 где черта  обозначает  замыкание оператора.}

Доказательство.  В утверждении 2 теоремы 5 положим $\psi(s)=c_1\cdot s$. Тогда по определению $\psi(A)=\overline{\sum_{j=1}^n c_1^j A_j}$. При этом мера  $\nu_t$ в представлении (\ref{gt})  есть $\varepsilon_{c_1^1t}\otimes\dots\otimes\varepsilon_{c_1^nt}$ ($\varepsilon_{b}$ --- мера Дирака, сосредоточенная в точке $b\in \mathbb{R}$), а потому в силу  (\ref{gt(A)})   полугруппа  $g_t(A)$ равняется  $\prod_{j=1}^n T_j( c_1^j t)$ и, следовательно,  непрерывна на
$(0;+\infty)$ в равномерной операторной топологии.

\textbf{Следствие 4.}  \textit{Если все полугруппы $T_j$  непрерывны на
$(0;+\infty)$ в равномерной операторной  топологии, а их генераторы $A_j$ не ограничены, то
\[
\sigma(A_j)={\rm pr}_j\sigma(A)={\rm pr}_j\sigma^{\prime\prime}(A), \  j=1,\dots, n,
\]
где ${\rm pr}_j$ --- ортогональная проекция пространства $\mathbb{C}^n$ на свое  $j$-е координатное подпространство. В частности, если спектр по крайней мере одного из операторов $A_j$ не пуст,   то и спектр $\sigma(A)$  не пуст.}

Для доказательства достаточно положить   в   предыдущем следствии   $c_1=e_j$,  где  $\{e_j, j=1,\dots, n\}$  --- стандартный базис в $\mathbb{R}^n$.

Для формулировки условия совпадения бикоммутантного спектра со спектром Шилова отождествим $\mathbb{C}^n$  с  $\mathbb{R}^{2n}$  с помощью отображения $z\mapsto(x_1,y_1,\dots,x_n,y_n)\ (z_k=x_k+iy_k)$ и рассмотрим конус  $K=\mathbb{R}^{2n}_+\subset\mathbb{C}^n$.  Напомним  (см., например, \cite[глава 1, \S 4]{KrNud}), что множество  $E\subseteq\mathbb{C}^n$ называется \textit{$K$-упорядоченно выпуклым}, если  $(b-K)\cap (a+K)\subseteq E$, как только   $a,b\in E$.

\textbf{Следствие 5.}  \textit{Если все полугруппы $T_j$  непрерывны на
$(0;+\infty)$ в равномерной операторной топологии,  их генераторы $A_j$ не ограничены,  а $\sigma(A)$  есть  $K$-упорядоченно выпуклое множество, то оно совпадает  с $\sigma^{\prime\prime}(A)$.}

Доказательство.  В силу теоремы 5.2 из главы 1 монографии \cite{KrNud},   множество $\sigma(A)$ есть пересечение двух семейств полупространств $(\lambda,z)\geq h_1(\lambda)\ (\lambda\in \Lambda)$  и $(\mu,z)\leq h_2(\mu)\ (\mu\in M)$, где $\Lambda$  и $M$ --- некоторые подмножества дуального конуса $K'$,  $h_1, h_2$ --- некоторые определенные на них функции, а $(\cdot , \cdot)$ --- стандартное скалярное произведение в $\mathbb{R}^{2n}$.  В нашем случае $K'=\mathbb{R}^{2n}_+$.   Пусть $z^0\in \sigma^{\prime\prime}(A)$. В силу следствия 3, для любого $\nu\in \mathbb{R}^{2n}_+$ найдется такая точка $z^1\in \sigma(A)$, что $\nu\cdot z^0=\nu\cdot z^1$.  Но тогда и $(\nu,z^0)=(\nu,z^1)$, а потому $(\lambda,z^0)\geq h_1(\lambda)\ (\lambda\in \Lambda)$  и $(\mu,z^0)\leq h_2(\mu)\ (\mu\in M)$, т.~ е. $z^0\in \sigma(A)$. Применение  теоремы 1  завершает доказательство.

\textbf{Следствие 6.}    \textit{Если все полугруппы $T_j$  непрерывны на
$(0;+\infty)$ в равномерной операторной топологии,  а их генераторы $A_j$ не ограничены, то при  $u\in \mathbb{R}_+^n,\ u\ne 0$ имеем}
\[
\sigma(T(u))= e^{u\cdot\sigma(A)}\cup\{0\}.
\]

Это вытекает из утверждения 3) теоремы 5 при $\psi(s)=e^{u \cdot  s} -1$.

Следующий результат была анонсирован, а также при дополнительных ограничениях доказан в \cite {SMZ} (см. там  теорему 12). Здесь мы снимаем эти ограничения.

 \textbf{Следствие 7.}  \textit{Для любой функции
$\psi\in {\cal T}_n$ и  набора $A$ попарно коммутирующих
генераторов ограниченных $C_0$-полугрупп в банаховом пространстве
$X$ справедливо включение}

\[
\sigma(\psi(A))\supseteq \psi(\sigma_J (A)).
\]

 Это  непосредственно вытекает из теоремы 1 и  теоремы 4.

\begin{center}
\textbf{5.   Приложение к исследованию асимптотического поведения многопараметрических полугрупп}
\end{center}

 В настоящее время  асимптотическое поведение однопараметрических полугрупп операторов довольно  хорошо изучено,  и соответствующая теория нашла важные применения  в теории уравнений в частных производных, теории вероятностей,  динамических системах, а также их приложениях  (см.  \cite{OP}, \cite{EngNag}, \cite{Em} и литературу там).  В данном разделе мы  применим полученные выше результаты к исследованию равномерной и сильной устойчивости  многопараметрических полугрупп.  Далее под \textit{конусом} в $\mathbb{R}^n_+$  мы будем понимать объединение лучей вида $l_u:=\{tu:t\in \mathbb{R}_+\}\ (u\in \mathbb{R}_+^n\setminus\{0\}$).

 \textbf{Определение  13.}  Будем говорить, что $n$-параметрическая полугруппа  $T$

 1) \textit{равномерно  экспоненциально устойчива  на конусе} $K$, если $\lim_{K\backepsilon v\to \infty}e^{\varepsilon|v|}\|T(v)\|=0$ для некоторого $\varepsilon>0$  (прямые скобки обозначают норму в $\mathbb{R}^n$);

 2) \textit{равномерно  устойчива  на конусе} $K$, если $\lim_{K\backepsilon v\to \infty}\|T(v)\|=0$;

 3) \textit{сильно экспоненциально устойчива  на конусе } $K$, если  для некоторого $\varepsilon>0$ $\lim_{K\backepsilon v\to \infty}e^{\varepsilon|v|}\|T(v)x\|=0$ при всех $x\in X$;

 4)  \textit{сильно устойчива  на конусе} $K$, если   $\lim_{K\backepsilon v\to \infty}\|T(v)x\|$ $=0$   при всех $x\in X$.

 В частности,  $T$  равномерно (сильно) устойчива  на луче $l_u$,  если $\|T(tu)\|\to 0$ при $t\to\infty$ (соответственно $\|T(tu)x\|\to 0$  при $t\to\infty$ для любого $x\in X$).

 В случае, когда $K=\mathbb{R}^n_+$,  вопрос об устойчивости
 полугруппы  $T$ на  $K$ тривиален --- устойчивость будет иметь место тогда и только тогда, когда каждая из частных полугрупп  $T_j$ устойчива на $\mathbb{R}_+$. Но, как показывает следующий пример, этот вопрос становится содержательным, если конус $K$ расположен строго внутри $\mathbb{R}^n_+$.

 \textbf{Пример 3.}  Пусть $X=C[0,1], q_1,q_2\in C[0,1]$, причем $q_1(0)=q_2(1)=0$,   $q_1(s)<0$  при $s\ne 0$, $q_2(s)<0$  при $s\ne 1$. Рассмотрим ограниченную 2-параметрическую полугруппу $T(t_1,t_2)x:=e^{t_1q_1+t_2q_2}x$ класса $(C_0)$ в $X$. Здесь частные полугруппы $T_j(t_j)x=e^{t_jq_j}x,\ j=1,2$ не являются сильно устойчивыми, поскольку, например, $|T_1(t_1)x(0)|=|x(0)|$. С другой стороны, полугруппа $T$  равномерно экспоненциально устойчива  на внутренних лучах конуса $\mathbb{R}^2_+$, так как для любого $k>0$ имеем $\|T(t,kt)\|=\max_{s\in[0,1]}|e^{t(q_1(s)+kq_2(s))}|=e^{t(q_1(s_0)+kq_2(s_0))}$ для некоторого $s_0\in [0,1]$.

 Для получения критериев равномерной экспоненциальной устойчивости многопараметрических полугрупп нам потребуется следующее понятие.

 \textbf{Определение  14.}  \textit{Спектральной гранью Шилова полугруппы} $T$ назовем число
\[
s(A)=\sup\left\{\sum\limits_{j=1}^n{\rm Re}\lambda_j:\lambda\in\sigma(A)\right\}.
\]
Заметим, что $s(A)\leq 0$, поскольку  $\sigma(A)\subset\times_{j=1}^n\sigma(A_j)$, а все полугруппы  $T_j$ предполагаются ограниченными.

В соответствии с известным обобщением классической  теоремы Ляпунова,  для однопараметрической $C_0$-полугруппы, становящейся  непрерывной  в равномерной операторной топологии,  равномерная устойчивость равносильна отрицательности ее спектральной грани  (см., например, \cite[c. 240]{OP}). Следующая теорема среди прочего содержит обобщение этого результата на случай многопараметрических полугрупп (ниже через  $r(B)$ обозначается спектральный радиус оператора $B$).

\textbf{Теорема 6.}  \textit{Пусть  все полугруппы $T_j$  непрерывны на
$(0;+\infty)$ в равномерной операторной топологии,  а их генераторы $A_j$ не ограничены. Следующие утверждения равносильны:}

(1)  $r(T(u))<1$ \textit{ для любого} $u\in (0,+\infty)^n$;

(2) $r(T(u))<1$ \textit{ для некоторого}  $u\in (0,+\infty)^n$;

(3) $s(A)<0$;

(4) \textit{полугруппа  $T$ равномерно устойчива на  луче  $l_u$  для любого}  $u\in (0,+\infty)^n$;

(5) \textit{полугруппа  $T$ равномерно устойчива на  луче  $l_u$  для некоторого}  $u\in (0,+\infty)^n$;

(6) \textit{полугруппа  $T$ равномерно устойчива на   любом замкнутом конусе}  $K\subset (0,+\infty)^n\cup\{0\}$;

(7) \textit{полугруппа  $T$ равномерно устойчива на некотором  замкнутом конусе}  $K\subset (0,+\infty)^n\cup\{0\}$.

(8) \textit{полугруппа  $T$ равномерно экспоненциально устойчива на   любом/некотором замкнутом конусе}  $K\subset (0,+\infty)^n\cup\{0\}$;

(9) \textit{полугруппа  $T$ сильно экспоненциально устойчива на   любом/некотором замкнутом конусе}  $K\subset (0,+\infty)^n\cup\{0\}$;

(10)  \textit{для любого/некоторого замкнутого конуса  $K\subset (0,+\infty)^n\cup\{0\}$ найдутся такие  числа}   $M_K>0$ \textit{ и} $\omega_K<0$, \textit{\textit{что}} $\|T(u)\|\leq M_Ke^{\omega_K|u|}$ \textit{при всех}  $u\in K$.

Доказательство будем проводить по схеме
\[
(1)\Rightarrow(2)\Rightarrow(3)\Rightarrow(1),\   (3)\Leftrightarrow(4)\Rightarrow(6)\Rightarrow(8)\Rightarrow(9)\Rightarrow(10)\Rightarrow(7)\Rightarrow(5)\Rightarrow(2).
\]
Утверждение (1)$\Rightarrow$(2) очевидно.   В силу следствия 6
\begin{equation}\label{logr}
\log r(T(u))=\sup\left\{\sum\limits_{j=1}^nu_j{\rm Re}\lambda_j:\lambda\in\sigma(A)\right\}.
\end{equation}
Следовательно, если   выполнено (2) и  $C=\max_j u_j$, то
\[
s(A)\leq\frac{1}{C}\sup\left\{\sum\limits_{j=1}^nu_j{\rm Re}\lambda_j:\lambda\in\sigma(A)\right\}<0.
\]

В свою очередь, если выполнено (3),  то при всех $u\in (0,+\infty)^n$ и $\lambda\in\sigma(A)$  имеем $\sum_{j=1}^nu_j{\rm Re}\lambda_j\leq ms(A)<0$,  где $m=\min_j u_j$, а потому $r(T(u))<1$ в силу (\ref{logr}).

Далее, условие (4) равносильно тому, что для любого   $u\in (0,+\infty)^n$ однопараметрическая полугруппа $U(t)=T(tu)$ с генератором $\overline{\sum_{j=1}^nu_jA_j}$  (см. доказательство следствия 3)  равномерно устойчива, что, ввиду  того, что она становится  непрерывной в равномерной операторной топологии,  равносильно отрицательности ее спектральной грани. С учетом следствия 3 заключаем, что условие (4) равносильно выполнению при всех   $u\in (0,+\infty)^n$ неравенства
\[
s\left(\overline{\sum\limits_{j=1}^nu_jA_j}\right)=\sup\left\{\sum\limits_{j=1}^nu_j{\rm Re}\lambda_j:\lambda\in\sigma(A)\right\}<0.
\]
Но выше  при доказательстве импликаций (2)$\Rightarrow$(3)$\Rightarrow$(1) было показано, что последнее условие  равносильно (3).

Покажем, что (4)$\Rightarrow$(6).   Допустим противное. Тогда существуют  замкнутый конус  $K\subset (0,+\infty)^n$ и число $\varepsilon>0$, такие,  что для любого натурального $N$  найдется вектор $v_N\in K$, для которого $|v_N|>N$ и $ \|T(v_N)\|\geq \varepsilon$. Заметим, что $v_N=t_Nu_N$, где   $u_N\in K_1$,  $K_1:=\{u\in K:|u|=1\}$ --- компакт,  и   $t_N>N$. Переходя, если нужно, к подпоследовательности, можем считать, что $u_N\to u_0\in K_1\ (N\to\infty)$. В силу (4),  найдется такое число $t_0>0$, что $\|T(t_0u_0)\|<\varepsilon/M^n$ (константа $M$  выбрана так, что $\|T_j(t)\|\leq M$  при всех $t\geq 0, j=1,\dots,n$). По причине  непрерывности   полугруппы $T$ на множестве $(0,+\infty)^n$ в равномерной операторной топологии,   для достаточно больших  $N$  выполняется неравенство  $\|T(t_0u_N)\|<\varepsilon/M^n$. Считая к тому же, что $N>t_0$, выводим  из этого неравенства, что
\[
\|T(v_N)\|=\|T(t_Nu_N)\|\leq\|T(t_0u_N)\|\|T((t_N-t_0)u_N)\|\leq\|T(t_0u_N)\|M^n<\varepsilon,
\]
а это  противоречит выбору $v_N$.

  Установим импликацию (6)$\Rightarrow$(8). Так как верно (4), то однопараметрическая полугруппа $t\mapsto T_1(t)\dots T_n(t)$ равномерно устойчива. А тогда, как известно (см., например, \cite[глава 5, предложение 1.3]{EngNag}), найдется такое
$\varepsilon>0$, что  $e^{\varepsilon t}\|T_1(t)\dots T_n(t)\|\to 0\ (t\to\infty)$. Следовательно, $n$-параметрическая полугруппа $u\mapsto e^{\varepsilon \sum_ju_j}T(u)$ равномерно устойчива на луче, а потому и на любом замкнутом конусе  $K\subset (0,+\infty)^n\cup\{0\}$ в силу (6).

Покажем теперь, что (9)$\Rightarrow$(10). Действительно, если полугруппа  $T$ сильно экспоненциально устойчива на    конусе  $K$, то  при некотором $\varepsilon>0$ для любого $x\in X$ множество $\{e^{\varepsilon |v|}\|T(v)x\|:v\in K\}$ ограничено, и по принципу равномерной ограниченности найдется такое число $M_K>0$, что $e^{\varepsilon|v|}\|T(v)\|\leq M_K$ при всех $v\in K$.

Импликации (10)$\Rightarrow$(7), (7) $\Rightarrow$ (5) и (8)$\Rightarrow$(9)   очевидны. Поскольку утверждение (5) $\Rightarrow$ (2)  следует из соответствующего одномерного результата (см., например, \cite[глава 5, предложение 1.7]{EngNag}),  теорема доказана.

Теорема 6 позволяет, в частности, применять к решению вопроса об  равномерной экспоненциальной устойчивости многопараметрической полугруппы критерии равномерной устойчивости однопараметрических полугрупп. В качестве примера приведем многомерную версию теоремы Ролевича (см. \cite[теорема 3.2.2]{Nee};  при $n=1, \varphi(t)=t^p, 1\leq p<\infty$  это известная теорема Датко-Паци).

\textbf{Теорема 7.} \textit{Пусть  все полугруппы $T_j$  непрерывны на
$(0;+\infty)$ в равномерной операторной топологии,  а их генераторы $A_j$ не ограничены.  Полугруппа $T$ будет  равномерно экспоненциально устойчива на   любом замкнутом конусе  $K\subset (0,+\infty)^n\cup\{0\}$, если найдется такой телесный замкнутый конус  $K_0\subset (0,+\infty)^n\cup\{0\}$ и такая неубывающая функция  $\varphi:\mathbb{R}_+\to\mathbb{R}_+$, что $\varphi(t)>0$ при $t>0$ и}
$$
\int\limits_{K_0}\varphi(\|T(v)x\|)dv<\infty\ \forall x\in X, \|x\|=1.
$$

Доказательство. Переходя, если нужно, к меньшему конусу, можем считать, что $K_0$ порождается  шаром $B[u_0,r]\subset \mathbb{R}^n\ (r>0)$.  Пусть $K_1=\{u\in K_0: |u|=1\}$  --- центральная проекция конуса $K_0$  на единичную сферу $S$ пространства   $\mathbb{R}^n$, $\mu$ --- сужение на  $K_1$  меры Лебега на  сфере $S$, $m$ --- линейная мера Лебега на  $\mathbb{R}_+$. Тогда $K_0=K_1\times \mathbb{R}_+$  и   $dv=\mu\otimes m$ --- cужение на $K_0$  $n$-мерной меры Лебега. Поэтому
$$
\int\limits_{K_0}\varphi(\|T(v)x\|)dv=\int\limits_{K_1}d\mu(u)\int\limits_{0}^\infty\varphi(\|T(su)x\|)ds\ (x\in X, \|x\|=1).
$$
Следовательно, для каждого $x\in X, \|x\|=1$ найдется такой     $u_x\in K_1$, что
 $$
\int\limits_{0}^\infty \varphi(\|T(su_x)x\|)ds<\infty.
$$
Пусть $u_x=t_xa_x$, где $t_x> 0, a_x\in B[u_0,r]$. Рассмотрим вектор $v_x=t_x(3u_0-a_x)$. Он принадлежит $K_0$, поскольку $|(3u_0-a_x)/2-u_0|<r$. Так как $u_x+v_x=3t_xu_0$ и $\|T(sv_x)\|\leq M^n$, то, полагая $\varphi_1(t):=\varphi(t/M^n)$, имеем
 $$
\int\limits_{0}^\infty \varphi_1(\|T(tu_0)x\|)dt=3t_x\int\limits_{0}^\infty \varphi_1(\|T(s3t_xu_0)x\|)ds=3t_x\int\limits_{0}^\infty \varphi_1(\|T(su_x+sv_x)x\|)ds=
$$
$$
= 3t_x\int\limits_{0}^\infty \varphi_1(\|T(sv_x)T(su_x)x\|)ds\leq 3t_x\int\limits_{0}^\infty \varphi(\|T(su_x)x\|)ds<\infty, \ \forall x\in X, \|x\|=1,
$$
а потому  $T$ равномерно устойчива на луче  $l_{u_0}$ по упомянутой теореме Ролевича. Применение теоремы 6 завершает доказательство.

 Что касается сильной устойчивости, то следующая простая лемма также позволяет сводить этот  вопрос к одномерному случаю.

\textbf{Лемма 4.} \textit{Пусть  $K$ --- замкнутый  конус в $(0,+\infty)^n\cup\{0\}$,  $x_0\in X$. Если  $\|T(tu)x_0\|\to 0\  (t\to\infty)$ для любого  $u\in K\setminus\{0\}$, то  $\|T(v)x_0\|\to 0\ (v\to\infty, v\in K)$.  В частности,  если полугруппа $T$ сильно устойчива на каждом луче $l_u, u\in K\setminus\{0\}$, то она сильно устойчива  на  $K$.}

Доказательство.   Пусть $v_m\in K, v_m\to\infty$. Тогда  $v_m=t_mu_m$, где   $u_m\in K_1, \ t_m\to+\infty\  (K_1=\{u\in K:|u|=1\})$. Как и в доказательстве теоремы 6  можно считать, что $u_m\to u_0\in K_1\ (m\to\infty)$.  По условию для любого $\varepsilon>0$ найдется такое число $t_0>0$, что $\|T(t_0u_0)x_0\|<\varepsilon/(2M^n)$. С учетом сильной непрерывности полугруппы $T$,  при достаточно больших $m$ имеем $\|T(t_0u_m)x_0-T(t_0u_0)x_0\|<\varepsilon/(2M^n)$. Если к тому же $m$ настолько велико, что $t_m>t_0$, то
\[
\|T(v_m)x_0\|=\|T((t_m-t_0)u_m)T(t_0u_m)x_0\|\leq M^n\|T(t_0u_m)x_0\|\leq
\]
\[
\leq M^n(\|T(t_0u_m)x_0-T(t_0u_0)x_0\|+\|T(t_0u_0)x_0\|)<\varepsilon,
\]
что и завершает доказательство леммы.

 Теперь легко получить для случая многопараметрических полугрупп версию известной теоремы Арендта-Бэтти-Любича-Ву  \cite{AB, LPh}  (см. также \cite[теорема 1.2.12]{Em}), дающей достаточные условия сильной   устойчивости ограниченных однопараметрических полугрупп.

\textbf{Теорема 8.} \textit{Пусть  $K$ --- замкнутый  конус в $(0,+\infty)^n\cup\{0\}$,  и пусть   все полугруппы $T_j$  непрерывны на $(0;+\infty)$ в равномерной операторной топологии,  а их генераторы $A_j$ не ограничены. Если для любого $u\in K\setminus\{0\}$  множество  $u\cdot\sigma(A)\cap i\mathbb{R}$ не более чем счетно, а $\sigma_R(\overline{u\cdot A})\cap i\mathbb{R}=\emptyset$, то полугруппа $T$ сильно устойчива  на  $K$.}

Доказательство. Как отмечено в доказательстве следствия 3,  однопараметрическая полугруппа
$$
U(t)=T(tu)=\prod\limits_{j=1}^n T_j(t u_j )\quad (u\in \mathbb{R}_+^n\setminus\{0\})
$$
  имеет генератор $G=\overline{u\cdot A}$. В силу того же следствия $\sigma(G)=u\cdot \sigma(A)$. Кроме того, точечный спектр сопряженного оператора есть $\sigma_p(G')=\sigma_R(\overline{u\cdot A})$. Следовательно, при любом $u\in K$ полугруппа $U$  сильно устойчива по упомянутой теореме Арендта-Бэтти-Любича-Ву,  и осталось применить лемму 4.

Аналогично, применяя следствие 5.1.13 из  \cite{Nee} (сильная устойчивость в \cite{Nee} называется равномерной) к полугруппам  $U$, рассмотренным выше, получаем следующее достаточное условие сильной устойчивости.

\textbf{Теорема 9.} \textit{Пусть  $K$ --- замкнутый  конус в $(0,+\infty)^n\cup\{0\}$,  и пусть   все полугруппы $T_j$  непрерывны на $(0;+\infty)$ в равномерной операторной топологии,  а их генераторы $A_j$ не ограничены. Если для любого $u\in K\setminus\{0\}$ выполнены следующие два условия:}

1) $\sigma_R(\overline{u\cdot A})\cap i\mathbb{R}=\emptyset$;

2) \textit{существует  плотное векторное подпространство  $Y_u\subset X$ со  свойством: для любого $y\in Y_u$ найдется такое счетное замкнутое $E\subset\mathbb{R}$, что   функция $\lambda\mapsto R(\lambda,\overline{u\cdot A})y$ допускает голоморфное продолжение в окрестность множества    $\{{\rm Re}\lambda\geq 0\}\setminus iE$,}

\textit{то полугруппа $T$ сильно устойчива  на  $K$.}

Таким же образом  с помощью \cite[следствие 5.1.16]{Nee} можно получить достаточные
 условия устойчивости индивидуальной орбиты. Полугруппа $T$  называется \textit{слабо почти периодической на  конусе К}, если  орбита $\{T(v)x: v\in K\}$ любого  элемента $x\in X$     относительно слабо компактна в $X$.

 \textbf{Теорема 10.} \textit{ Пусть  $K$ --- замкнутый  конус в $(0,+\infty)^n\cup\{0\},\ x_0\in X$,  и    все полугруппы $T_j$  непрерывны на $(0;+\infty)$ в равномерной операторной топологии,  а их генераторы $A_j$ не ограничены. Если полугруппа  $T$ слабо почти периодична на К, для любого $u\in K\setminus\{0\}$ справедливо равенство $\sigma_p(\overline{u\cdot A})\cap i\mathbb{R}=\emptyset$,  и найдется такое счетное замкнутое $E_u\subset\mathbb{R}$, что   функция $\lambda\mapsto R(\lambda,\overline{u\cdot A})x_0$ допускает голоморфное продолжение в окрестность множества    $\{{\rm Re}\lambda\geq 0\}\setminus iE_u$, то $\|T(v)x_0\|\to 0$ при  $v\in K, v\to\infty$.}

 Следующий пример показывает, как предыдущие теоремы  могут быть применены в ситуации, которую можно назвать "каскадом задач Коши". Для простоты изложения мы рассматриваем лишь случай $n=2$.

\textbf{Пример 4.}  Пусть $A_1, A_2$ --- коммутирующие неограниченные генераторы ограниченных $C_0$-полугрупп $T_1, T_2$ соответственно в банаховом пространстве  $X$. Предположим,  что $D(A_2)\subseteq D(A_1)$ и что  $T_1$ и $T_2$    непрерывны на $(0,\infty)$ в равномерной операторной  топологии (например, голоморфны), и рассмотрим задачу Коши
$$
 \left\{\matrix {\frac{\partial u(t_1,t_2)}{\partial t_1}=A_1u(t_1,t_2), &
(t_1,t_2)\in \mathbb{R}_+^2, \cr
\hspace{-10mm} u(0,t_2)=v(t_2) \cr}\right.\eqno(18)
$$
с неизвестной функцией $u: \mathbb{R}_+^2\to X$.  Пусть функция $v:\mathbb{R}_+\to X$, в свою очередь, является решением задачи Коши

$$
 \left\{\matrix {\frac{d v(t_2)}{d t_2}=A_2v(t_2), &
t_2\in \mathbb{R}_+, \cr
\hspace{-7mm} v(0)=v_0, & v_0\in D(A_2). \cr}\right.\eqno(19)
$$
Как известно, решение задачи (19) дается формулой $v(t_2)=T_2(t_2)v_0$, причем $v(t_2)\in D(A_2)\subseteq D(A_1)$ при всех $t_2\in \mathbb{R}_+$. Следовательно, решение задачи (18) имеет вид $u(t_1,t_2)=T(t_1,t_2)v_0$, где  $T(t_1,t_2)=T_1(t_1)T_2(t_2)$.

Пусть $K$ --- (замкнутый) угол с вершиной в начале координат, содержащийся в $(0,\infty)^2\cup\{0\}$. Согласно теореме 6, если (и только если) $s(A_1,A_2)<0$ (или $r(T(t_1^0,t_2^0))<1$  для некоторой точки  $(t_1^0,t_2^0)\in (0,+\infty)^2$), то   найдутся такие  числа   $M_K>0$  и $\omega_K<0$, что  $\|u(t_1,t_2)\|\leq M_Ke^{\omega_K(t_1+t_2)}\|v_0\|$ при всех $(t_1,t_2)\in K$. То же заключение справедливо при выполнении условий теоремы 7.

Если же  для любого $(t_1,t_2)\in K$  множество  $(t_1,t_2)\cdot\sigma(A_1,A_2)\cap i\mathbb{R}$ не более чем счетно, а $\sigma_R(\overline{t_1A_1+t_2A_2})\cap i\mathbb{R}=\emptyset$, то $\|u(t_1,t_2)\|\to 0$   для любого $v_0\in D(A_2)$ при $(t_1,t_2)\to\infty,\ (t_1,t_2)\in K$ по теореме 8.  То же заключение справедливо при выполнении условий теоремы 9.

\textbf{Благодарности.} Автор признателен рецензентам за замечания и предложения, учет которых привел к существенному улучшению рукописи. Автор благодарит также А.~Г.~Баскакова,  любезно предоставившего ему оттиски своих работ \cite{Bas1} и \cite{Bas2}.

\newpage

\centerline{СВЕДЕНИЯ ОБ АВТОРЕ}

Миротин Адольф Рувимович,

полное название научного учреждения: Гомельский гос.  университет им. Ф. Скорины,

ул. Советская, 104 ,  г. Гомель
246019, Беларусь.

Адрес для переписки: 246028, Беларусь, г. Гомель, ул.
Старочерниговская, д. 6, кв. 65

Тел.: +(375232)579511 (дом.),  +375293988656 (моб.),

+(375232)573051 (служ.)

 Адрес электронной почты:
amirotin@yandex.ru

\begin{center}
{\bf A. R. Mirotin \\ ON    JOINT  SPECTRA  OF FAMILIES OF UNBOUNDED OPERATORS  }
\end{center}

\end{document}